\documentstyle[12pt,epsf]{article}

\def\reel{\hbox{{\rm R}\kern-1em\hbox{{\rm I} }}}
\def\relatif{\ \hbox{{\rm Z}\kern-.4em\hbox{\rm Z}}}
\def\nat{\hbox{{\rm N}\kern-1em\hbox{{\rm I} } }}
\def\comp{\hbox{{\rm C}\kern-.55em\hbox{{\rm I} } }}
\def\smallcomp{\hbox{\fiverm C}\kern-.35em{\hbox{\fiverm I}}}

\def\fudge{\mathchoice{}{}{\mkern.5mu}{\mkern.8mu}}
\def\bbc#1#2{{\rm \mkern#2mu\vbar\mkern-#2mu#1}}
\def\bbb#1{{\rm I\mkern-3.5mu #1}} \def\bba#1#2{{\rm #1\mkern-#2mu\fudge
#1}}
\def\bb#1{{\count4=`#1 \advance\count4by-64 \ifcase\count4\or\bba
A{11.5}\or
\bbb B\or\bbc C{5}\or\bbb D\or\bbb E\or\bbb F \or\bbc G{5}\or\bbb H\or
\bbb I\or\bbc J{3}\or\bbb K\or\bbb L \or\bbb M\or\bbb N\or\bbc O{5} \or
\bbb P\or\bbc Q{5}\orrrr\b
bb R\or\bbc S{4.2}\or\bba T{10.5}\or\bbc U{5}\or    \bba V{12}\or\bba
W{16.5}\or\bba X{11}\or\bba Y{11.7}\or\bba Z{7.5}\fi}}
\def\rat{\hbox{{\rm Q}\kern-.70em\hbox{{\rm I} } }}

\def \Z {\relatif}

\newcounter{theoremcounter}
\newcounter{lemmacounter}

\newtheorem{thm}[theoremcounter]{Theorem}
\newtheorem{lemma}[lemmacounter]{Lemma}

\newcommand{\be}{\begin{equation}}
\newcommand{\ee}{\end{equation}}

\newcommand{\ber}{\begin{eqnarray}}
\newcommand{\eer}{\end{eqnarray}}

\newcommand{\nin}{\noindent}
\newcommand{\non}{\nonumber}

\def\qed{\hfill \vrule height1.3ex width1.2ex depth-0.1ex}
\def\bbb#1{{\rm I\mkern-3.5mu #1}} \def\bba#1#2{{\rm #1\mkern-#2mu\fudge
#1}}

\title{Asymptotic formula for a partition function of reversible
coagulation -fragmentation processes.}

\author{{\bf Gregory A.Freiman}
\thanks{E-mail: grisha@math.tau.ac.il }\\
School of Mathematical Sciences, Raymond and Beverly Sackler
Faculty\\ of Exact Sciences, Tel-Aviv University,Ramat-Aviv,\\
Tel-Aviv, Israel.\\ \quad {\bf Boris L. Granovsky}
\thanks{E-mail: mar18aa@techunix.technion.ac.il} \\
Department of Mathematics, Technion-Israel Institute of Technology,\\
Haifa, 32000,Israel.
}
\begin{document}
\maketitle
\vskip 5cm

\nin American Mathematical Society 1991 subject classifications.

\nin Primary-60J27;secondary-60K35,82C22,82C26.

\nin Keywords and phrases: Coagulation-Fragmentation process,
Local limit theorem,
Distributions on the set of partitions.

\newpage
\setcounter{equation}{0}
\begin{abstract}
We construct a probability model seemingly unrelated to the
considered stochastic process of coagulation and fragmentation. By
proving for this model the local limit theorem, we establish the
asymptotic formula  for the partition function of the equilibrium
measure  for a wide class of parameter functions of the process.
This formula proves the conjecture stated in \cite{dgg} for the
above class of processes. The method used goes back to
A.Khintchine.
\end{abstract}
\section{Introduction and Summary.}

\nin The motivation for our  research came from a conjecture
stated in \cite{dgg}, p.462, in  the following setting.

 \nin For a
given integer $N,$ denote by
\be
\eta=( n_1, \ldots, n_N): 0\le n_k\le N,
\quad \sum_{k=1}^{N} k n_k = N
\label{0}
\ee
a partition of  $N$
into $n_k$ groups of size  $k,
\quad k=1, 2, \ldots,N$ and by $\Omega_N = \{\eta \}$ the set of all
partitions of $N.$

\nin We will be interested in the particular probability measure
$\mu$ on $\Omega_N$ given by

\ber \mu_N (\eta) = C_N \frac{a_1^{n_{1}} a_2^{n_{2}} \ldots
a_N^{n_{N}}} {n_{1}! n_{2}! \ldots n_{N}!}, \quad
\eta=(n_1,\ldots,n_N)\in \Omega_N, \non \\ \label{mu} \eer

\nin where $a_k>0, $  $k=1,2,\ldots, N$  and $C_N =
C_N(a_1,\ldots,a_N)$ is the partition function of the distribution
$\mu_N:$

\ber c_N: = C_N^{-1} = \sum_{\eta\in \Omega_N} \frac{a_1^{n_1}
a_2^{n_2} \ldots a_N^{n_N}} {n_1!n_2!  \ldots n_N!}, \quad
\eta=(n_1,\ldots,n_N)\in \Omega_N, \nonumber\\
 \quad N = 1, 2, \ldots
\label{1} \eer

\nin The probability measure $\mu_N$ is  the equilibrium state of
a class of reversible coagulation-fragmentation processes
(CFP's)(see for references \cite{dgg}).

\nin CFP's  trace their history from Smoluchowski(1916) and they
have been intensively studied since  this date. The process models
the stochastic evolution in time of a population of $N$ particles
distributed into groups that coagulate and fragment at different
rates. The model arises in different contexts  of
application:polymer kinetics, astrophysics, aerosols,  biological
phenomena, such as animal grouping, blood cell aggregation, etc.
Observe that particular choices of $a_k, \  k=1,2, \ldots, N$ in
(\ref{mu}) lead to a variety of known stochastic models. For
example, when $a_k=\frac{\beta}{k}, \  k=1,2, \ldots, N, \
\beta>0,$ (\ref{mu}) becomes the widely known Ewens sampling
formula that arises in population genetics.

\nin  Following \cite{dgg}, we view CFP as a continuous-time
Markov process on the state space $\Omega_N.$ Formally, a CFP is
given by the rates $\psi$ and
 $\phi$ of the two
 possible transitions:
coagulation and fragmentation respectively. Namely, $\psi(i,j),
\quad 2\le i+j \le N$ is the rate of merging  of two groups of
sizes $i$ and $j$ into one group of size $i+j, $ and
$\phi(i,j),\quad 2\le i+j \le N$ is the  rate of splitting of a
group of size $i+j$ into two groups of sizes $i$ and $j.$  We
consider the class of CFP's for which the ratio of the transition
rates has the form

\ber  \frac{\psi(i,j)}{\phi(i,j)}=\frac{a_{i+j}}{a_ia_j}, \quad
i,j: 2\le i+j \le N, \non
\\ \label{q} \eer

\nin where $a_k>0,  k=1,\ldots,N$  are given parameters of the
process. Owing to (\ref{q}), the condition of detailed balance
holds, and, consequently, the CFP considered is reversible with
respect to the  invariant measure (\ref{mu}).

\nin Letting $N\to \infty,$  we will be concerned with the
relationship  between two infinite sequences $\{a_n\}_1^\infty$
and $\{c_n\}_0^\infty, c_0=1.$

\nin It was conjectured in \cite{dgg}* that the existence of the
limit

\be
 \lim_{n \to \infty} \frac{a_n}{a_{n+1}} >0
\label{rs}
\ee

\nin implies the existence of the limit
\be \lim_{n \to \infty} \frac{c_n}{c_{n+1}} >0.
\label{rs1}
\ee

\nin Apart from the fact that the conjecture is a challenging
mathematical problem, one can see from \cite{dgg} that it also has
a direct significance for the stochastic model in question. First,
if the limit (\ref{rs1}) exists, then a variety of functionals of
the process( e.g., the expected values and variances of finite
group sizes),  as $N\to \infty$, can be explicitly expressed via
this limit. Next, by formula (4.16) in \cite{dgg} we have that
\ber cov(n_k, n_l) = a(k)a(l) \left(\frac {c_{N - k -l}}{c_N} -
\frac{c_{N-k} c_{N-l}}{(c_N)^2}\right), \nonumber\\ k\neq l =
1,2,\ldots, N,\quad k+l\le N. \label{rnkl} \eer

\nin Thus, the validity of the conjecture will imply that at the
steady state  the random variables $n_k, n_l, \ \ k\ne l$ become
uncorrelated, as $N\to\infty.$ This fact incorporates into the
assumption of independence of  sites in  mean-field models, as
$N\to \infty,$ that is commonly accepted in statistical
 physics.

\nin Another motivation for our study is provided by a quite
different  field, known as  random combinatorial
structures(RCS's). The connection of CFP's to this field  is based
on the following observation made in \cite{dgg}. Let  $Z_i,
i=1,\ldots,N$ be independent Poisson random variables with
respective means $a_i>0,\quad i=1,\ldots,N.$ Then it is easy to
see that the distribution $\mu_N$ admits the following
representation

\ber \mu_N(\eta)= Pr\{ Z_1 =n_1, \ldots, Z_N=n_N \vert
\sum_{i=1}^{N} i Z_i = N \},\nonumber\\
\eta=(n_1,\ldots,n_N)\in\Omega_N. \label{rep} \eer \vskip .5cm

\nin It turns out that (\ref{rep}) is the general form of
distributions arising in a variety of RCS's.  This is explained in
\cite{arr}, \cite{arr2} and \cite{kol}. (Theorem 1, p.96 in
\cite{arr} gives a rigorous proof of this fact). The simplest
example of a RCS is  a random choice from $N!$ permutations of $N$
objects. Cauchy's formula for the number of permutations  having
$n_k$ cycles of length $k, \ k=1,\ldots,N,$ where $\sum_{k=1}^{N}
k n_k = N,$ tells us that the probability of picking a permutation
with this property  is given by (\ref{mu}) with $a_k=k^{-1}, \
k=1,\ldots, N$ and $C_N=1.$

Added in proofs:The conjecture was recently proved by J.Ball and
S.Burris in" Asymptotics for Logical Limit Laws", 2001, Preprint.

\nin In view of this, one can translate the preceding reasoning in
the context of the Conjecture, into the language of RCS's.

 \nin Our paper is devoted exclusively to the study of the
asymptotic behaviour, as $n\to \infty,$ of the quantity $c_n,$
defined by (\ref{1}).

\nin The asymptotic formula for $c_n$ established in our paper
proves the Conjecture for a wide class of parameter functions $a:
a(k)=a_k, \quad k=1,2,\ldots$. We mention also two other
applications  of our result related to global characteristics of
CFP's(=RCS's).

\nin (i) For a given $n,$ denote by $v_n$ the mean value of the
total number of different groups at the equilibrium of CFP
(=components in a RCS). It follows from (4.15) in \cite{dgg} that
\be
v_n=\sum_{k=1}^n a_k \frac{c_{n-k}}{c_n}. \label{1vn}\ee

\nin (ii) Denote by $p{_\infty}$ the probability at the steady
state of the creation of a cluster
 (=  component in a RCS) of infinite size. It was shown in
\cite{dgg}, p.462, that the condition \be \lim_{n\to \infty}( v_n-
v_{[\alpha n]-1})>0, \label{vn}\ee for some $0<\alpha\le 1$ is
sufficient for $p_{\infty}>0.$ (Here $[\bullet]$ is the integer
part of a number.)

\nin Thus, with the help of the asymptotic formula for $c_n,$ one
can reveal the asymptotic behaviour of $v_n,$ as $n\to \infty$
and, consequently,  find the limit in (\ref{vn}). The latter will
answer a question which is common in statistical physics.

 \nin Also we want to point out that  determining the asymptotic
 properties of  partition
functions for interacting particle systems is a difficult
mathematical problem widely discussed  in statistical physics (see
e.g. \cite {th}).

\section{Description of the method and\\
 a sketch of its history}

\nin We assume  $ a_n >0,  n=1,2,\ldots$ and  that the following
limit  exists:
\be
 \lim_{n \to \infty} \frac{a_n}{a_{n+1}}:=R, \quad 0<R<\infty.
\ee

\nin Thus, the power series in $x,$

\be
S(x): = \sum_{n=1}^{\infty} {a_n x^n}, \label{gen1} \ee

\nin has  radius of convergence $R$ and it converges in the
complex domain $D\subseteq \{\vert x\vert \le R\}.$

\nin Then (see \cite{dgg}), $g(x)=e^{S(x)}, \   x\in D$ is the
generating function for the  sequence $\{c_n\}_0^\infty$ defined
by
 (\ref{1}).
Namely,
\be
g(x)= e^{S(x)} = \sum_{n=0}^{\infty} c_n x^n, \quad x\in D,
\label{gen2} \ee

\nin and moreover, the series (\ref{gen1}) and (\ref{gen2})
converge in the same domain $D$.

 \nin The method we use here for
deriving the asymptotic formula for $c_n$ goes back to A.
Khintchine's pioneering monograph \cite{Kh}. In \cite{Kh}
Khintchine developed the idea of expressing of
 values of quantum

 \nin statistics via  the probability function of a sum
 of correspondingly constructed independent integer- valued
 random variables.
Subsequent implementation of the local limit theorem resulted in
the method of the derivation of asymptotic
 distributions of quantum
statistics. In \cite{Kh}
  this method was systematically
applied to systems of photons and some other models. The method
was further developed by A. Postnikov and G. Freiman, (see for
references \cite{po}) who applied it to analytic number theory. In
particular, G. Freiman formulated a local limit theorem for some
 asymptotic problems related to partitions. A general scheme for the
 derivation
 of asymptotic formulae for these kind of problems was outlined
 by G.Freiman and J. Pitman in
 \cite{fr},\cite{fr1} ( for references see also  \cite{as}.)

\nin A similar approach, also based on the implementation of the
local limit theorem, has been independently developed for the last
fifteen years in the theory of RCS's. A very good exposition of
this direction of research is given in the  recent monograph
\cite{kol} by V. Kolchin. We will explain briefly the  basic
difference between the problem addressed in the present paper and
those  in \cite{kol}. In the context of the generalized scheme of
allocation that encompasses a variety  of RCS's, $S$ and $g$ are
the generating functions for, respectively, the total number of
combinatorial objects of size $n$ and for the number of such
objects possessing  a definite property.  In this setting it  is
assumed that the expression for the function $g$ is known
explicitly. Based on this, the combinatorial quantity in question
is expressed via the probability function of a sum of i.i.d.
discrete random variables, distributed according to a probability
law that depends
 on the given values of  $\{c_n\}_0^\infty.$ Such scheme is applicable
for example, for investigation  of the asymptotic of the number
$F_{n,N}$ of all forests of $N$ nonrooted trees having $n$
vertices, in which  case $c_n=(n!)^{-1}n^{n-2}, n=1,2,\ldots.$

\nin First, w.l.o.g. we assume throughout the paper  that the
common radius of convergence
 of the series (\ref{gen1}) and (\ref{gen2}) equals  $ 1.$
This is due to the fact that taking in (\ref{1}) $\tilde a_j= R^j
a_j,\quad j=1,\ldots, N,$ it follows from (\ref{0})  that $\tilde
c_N= R^N c_N.$ The above assumption makes the $\lim_{n\to
\infty}\frac{a_n}{a_{n+1}},$ if it exists, equal  $1.$

\nin Our starting point is the following representation of $c_n.$

\begin{lemma} \hfill
\label{lem:genfunction}

\nin
\be
c_n= e^{n\sigma}\int _0^1  \prod_{l=1}^{n}\left(
\sum_{k=0}^{\infty}\frac{a_l^k
e^{-lk\sigma+2\pi i\alpha lk}}{k!}\right)\times e^{-2\pi
i\alpha n}d\alpha, \quad n=1,2,\ldots
\label{cn}
\ee
\nin for any  real $\sigma.$
\end{lemma}

\bigskip
\nin
{\bf Proof:}

\nin \nin It follows from (\ref{1}) that $c_n$ depends only on
$a_1,\ldots, a_n,$ which means that the first $n+1$ terms of the
Taylor series expansions of the two functions

\be
g(x)=e^{S(x)} \quad {\rm and} \quad   g_n(x):= e^{\sum_{l=1} ^{n}
a_lx^l}, \quad x\in D \ee are the same, i.e. $c_k=c_{k,n}, \quad
k=0,\ldots,n$, where $\{c_{k,n}\}_{k=0}^\infty$ is the sequence
related to the function $g_n.$ \nin For a fixed $n$, the series
expansion of the function $g_n(x)$ converges for all $x.$ So,  we
can set

\be
x=e^{-\sigma +2\pi i\alpha},
\label{-1}
\ee
for some real $\sigma $ and   $\alpha.$

\nin Then we  have

\be
\int _0^1 g_n(x) e^{-2\pi i \alpha n} d\alpha = \int _0^1 \left
(\sum_ {k=0}^{\infty} c_{k, n} e^{-\sigma k+2\pi i \alpha
(k-n)}\right ) d\alpha =c_n e^{-n\sigma}. \label{2} \ee \nin
 The last equality is due to the fact that

\be
\int _0^1 e^{2\pi i \alpha m}d\alpha =
\left \{\begin{array}{ll}1, & {\rm if~} m=0 \cr
0,& {\rm if~} m\neq 0,\quad m\in \Z.
\end{array}
\right.
\ee

\nin Finally, substituting
\be
g_n(x)= \prod_{l=1}^{n} e^{a_lx^l}= \prod_{l=1}^{n}
\left(\sum_{k=0}^{\infty}
\frac{(a_lx^l)^k}{k!}\right)
\ee
and (\ref{-1}) in the LHS of (\ref{2}) we get the claim. \qed

\nin Our next step will be to give a probabilistic meaning to the
expression (\ref{cn}) for $c_n.$

\nin We introduce the following notations.

\be
S_n(x)=\sum_{l=1} ^{n} a_lx^l,
\ee

\be p_{lk}=\frac {a_{l}^k e^{-\sigma lk}}{k!\exp\left(
a_le^{-\sigma l}\right)}, \quad l=1,\ldots,n,\quad k=0,1,\ldots
\label{p} \ee

\be
\varphi_l(\alpha)=  \sum_{k=0}^{\infty} p_{lk}e^{2\pi i\alpha lk},
\quad \alpha\in R, \label{v1} \ee

\be
\varphi(\alpha)= \prod_{l=1}^n \varphi_l(\alpha) , \quad \alpha\in
R . \label{v2} \ee

\nin Now (\ref{cn})  can be rewritten as
\be
c_n= e^{n\sigma} e^{S_n(e^{-\sigma})}\int_0^1
\varphi(\alpha)e^{-2\pi i\alpha n}d\alpha.
\label{3}
\ee

\nin The fact that for a given $l(1\le l\le n ),$ $p_{lk}, k=0,1,
\ldots $ is  a Poisson probability function with parameter
$a_le^{-\sigma l},$  suggests  the following probabilistic
interpretation of the integral in the RHS of (\ref{3}).

\nin Let $X_1,\ldots,X_n$ be  independent integer-valued random
variables defined by

\be
Pr(X_l=lk)=p_{lk}, \quad l=1,\ldots, n , \quad k=0,1,\ldots
\label{pr} \ee

\nin Then $\varphi(\alpha)$ defined above is the characteristic
function of the sum
$Y=X_1+\ldots+X_n$ and we have

\be
\int_{0}^{1} \varphi(\alpha)e^{-2\pi i\alpha n}d\alpha =Pr(Y=n).
\label{333} \ee \vskip .5cm \nin Now (\ref{3}) can be viewed as an
analog of the aforementioned Khintchine's representation for
$c_n.$

\nin  It is well- known \cite{gn}
  from the classical theory of limit distributions
of sums  of independent  integer-valued
random variables
 that,
 under certain conditions  on distributions of the variables,
a  local limit theorem is valid.

\nin  In our subsequent study, the free parameter $\sigma$ will be
taken depending on $n: \sigma=\sigma_n.$  By (\ref{pr}) and
(\ref{p}) this means that the probability law of each of the $n$
random variables $X_l, \quad l=1,2,\ldots, n$ depends on $n.$
Therefore, to compare with the classical case, we  will be dealing
here with a triangular array of random variables. For this case,
general necessary and sufficient conditions  for validity of the
local limit theorem are not known. For some  cases
  results in this direction
were obtained in \cite{ms} and \cite{des}. In the first of these
two papers a sufficient condition  was established (see
(\cite{ms}, Theorem 2, condition III) in the case of  an array  of
general lattice random variables. It can be verified that this
condition ( which can be viewed as a version of the celebrated
condition of asymptotic uniformity)
  fails for the class of parameter functions $a$
considered in our paper. The second paper studies  exclusively
 the case of the triangular array of trinomial random variables.

\nin Most of this paper is devoted to the proof of the local limit
theorem in  the  above setting.

\nin Namely, we will demonstrate that under certain conditions on the
parameter function $a$

\be
Pr(Y=n)\sim (2\pi B_n^2)^{-1/2}e^{-(M_n-n)^2/2B_n^2}, \quad {\rm
as} \quad n\to \infty, \label{4} \ee

\nin where $M_n=EY$ and $B_n^2=Var Y.$

\section{Proof of the local limit theorem}

\nin  In order to prove (\ref{4}) we have to find the asymptotic
formula, as $n\to \infty,$ for the integral in the RHS of
(\ref{3}). We will denote in the sequel, by $\gamma, \gamma_i, \
i=1,2,\ldots$ constants.

 \nin First, we obtain  the explicit expressions for the
quantities $M_n$ and $B_n^2.$

\nin  (\ref{pr}) and (\ref{p}) say that $l^{-1}X_l, \l=1,\ldots,$
are Poisson($a_le^{-\sigma l}$) random variables. So,  we have

\be
EX_l= la_l e^{-\sigma l},\quad l=1,\ldots,n , \label{abc} \ee

\be
Var X_l=l^2a_le^{-l\sigma}, \quad l=1,\ldots,n. \ee

\nin This gives

\begin{eqnarray}
M_n &=&\sum_{l=1}^n la_l e^{-l\sigma}, \quad n=1,2,\ldots \\
B_n^2 &=& \sum_{l=1}^n l^2a_l e^{-l\sigma }, \quad n=1,2,\ldots
\end{eqnarray}

\nin It follows from the preceding discussion that the
representation (\ref{3}) holds for any real $\sigma.$ Our next
result shows that  $\sigma$ can be chosen so  that the exponential
factor  in the RHS of (\ref{4}) equals  1, for any $n=1,2,\ldots$.

\begin{lemma} \hfill
\label{lm1}

\nin The equation

\be
\sum_{l=1}^n la_l e^{-l\sigma}=n
\label{sig}
\ee

\nin has a unique solution $\sigma=\sigma_n,$ for any $n=1,2,\ldots.$
\end{lemma}
\nin {\bf Proof:} The assertion follows immediately from the
assumption $a_l> 0, \quad l=1,2,\ldots $ \qed

\nin {\bf Remark} The above choice of the free parameter $\sigma$
makes  the probability of the event $\{Y=n\}$ large, as $n\to
\infty$. The same idea is widely used  for approximation of RCS's
by independent processes (see for references \cite{arr}). In
statistical physics such a way of choosing a free parameter for
estimating  averages is known as Darwin-Fawler method developed in
the $1930's$  (for references see, e. g. \cite{stat})

\nin It follows from (\ref{sig}) that if  the series
$\sum_{l=1}^{\infty} la_l$ converges, then $\sigma_n \le 0 $ for
sufficiently large $n$, while in the opposite case the sign of
$\sigma_n$ depends on the behaviour of $S_n^{'}(1),$ as $n\to
\infty.$ However, in both cases the following basic property of
$\sigma_n, \quad n=1,2,\ldots$  holds.

\begin{lemma} \hfill
\label{lm2}

\nin Let
\be
 \lim_{n \to \infty} \frac{a_n}{a_{n+1}}=1.
\label{5}
\ee

\nin Then $\lim_{n\to \infty} \sigma_n =0.$
\end{lemma}
\nin
{\bf Proof:}

\nin By the definition of $\sigma_n,$

\be
n^{-1}\sum_{l=1}^n la_l e^{-l\sigma_n}=1, \quad n=1,2,\ldots
\label{7}
\ee

\nin Denote $b_l=la_l>0, \ l=1,2\dots.$ Based on (\ref{5})

\nin let $N=N(\epsilon), \quad \epsilon >0 $ be s.t.

\be
1-\epsilon\le \frac{b_{l+1}}{b_l}\le 1+\epsilon, \quad \mbox{ for
all  \quad $l\ge N.$} \ee

\nin Consequently,
\be
(1-\epsilon)^{l-N}b_N\le b_l\le (1+\epsilon)^{l-N} b_N, \quad l\ge N.
\ee

\nin Now suppose that $\lim_{k\to \infty} \sigma_{n_k}=\sigma,$
for some subsequence $ n_k\to +\infty,$ as $k\to \infty,$ where
$\vert\sigma\vert\le\infty.$

\nin Let first $-\infty\le\sigma
<0,$ then taking $0<\epsilon
< 1-e^{\sigma /2}$ we have
$$ n_k^{-1}\sum_{l=1}^{n_k} b_l e^{-l\sigma_{n_k}}\ge
(1-\epsilon)^{-N}b_N
n_k^{-1}\sum_{l=N}^{n_k}[(1-\epsilon)e^{-\sigma_{n_k}}]^l\to\infty,
$$ \be \mbox{as $k\to \infty$}, \label{222} \ee \nin since
$$(1-\epsilon)e^{-\sigma_{n_k}}>
e^{\frac{\sigma}{2}-\sigma_{n_k}}\to e^{-\sigma/2}>1,\quad \rm
{as} \quad k\to \infty.$$

\nin If now $0<\sigma \le \infty,$ then for $0<\epsilon< e^{\sigma
/2} -1$,

\begin{eqnarray}
 \lefteqn{n_k^{-1}\sum_{l=1}^{n_k} b_l e^{-l\sigma_{n_k}}\le {}}
\nonumber\\& &{{ }n_k^{-1}\left(\sum_{l=1}^{N} b_l
e^{-l\sigma_{n_k}}+(1+\epsilon)^{-N}b_N \sum_{l=N}^{n_k}\left(
(1+\epsilon) e^{-\sigma_{n_k}}\right)^{l}\right )\to 0,{ }
}\nonumber\\ & &\mbox{as \quad $k\to \infty$}, \label{111}
\end{eqnarray} \nin since in this case  $e^{-\sigma/2}<1.$ Both
(\ref{222}) and (\ref{111}) contradict (\ref{7}), which implies
that  $\limsup_{n\to \infty} \sigma_n=\liminf_{n\to
\infty}\sigma_n=0.$ \qed \vskip .5cm

\nin In what follows we will assume that $\sigma =\sigma_n,$
as defined by (\ref{sig}).
\nin Our next lemma provides the expression for
the integrand in the LHS of (\ref{333})
for small values of $\alpha$.

\begin{lemma} \hfill
\label{lm3}

\nin For a fixed $n$ and  $\sigma=\sigma_n,$

\be
\varphi(\alpha)e^{-2\pi i\alpha n} =
\exp\left(-2\pi^2\alpha^2B_n^2 +
 O(\alpha^3\rho_3)\right), \quad
\mbox{as $\alpha\to 0$}, \label{10} \ee \nin where
$\rho_3=\rho_{3}(n)=\sum_{l=1}^n l^3a_l e^{-l\sigma_n },\quad
n=1,2,\ldots.$
\end{lemma}

\nin {\bf Proof:} By (\ref{p}) - (\ref{v2}),

$$ \varphi_l(\alpha)= \sum_{k=0}^{\infty} \frac{\left(a_l
e^{-\sigma_nl}e^{2\pi i\alpha l}\right )^k}
{k!\exp\left(a_le^{-l\sigma_n}\right)}=
\exp\left(a_le^{-l\sigma_n}\left(e^{2\pi i\alpha l}-1\right
)\right),$$ \be \alpha\in R \label{8} \ee

\nin  and
\be
\varphi (\alpha)= \exp\left(\sum_{l=1}^n
a_le^{-l\sigma_n}\left(e^{2\pi i\alpha l} -1\right )\right),
\alpha\in R. \label{9} \ee

\nin Finally, substituting  in  (\ref{9})  the Taylor expansion
(in $\alpha$)

\be
e^{2\pi i\alpha l}-1= 2\pi i\alpha l-2\pi^2\alpha^2l^2
+O(\alpha^3l^3), \quad \mbox{as $\alpha\to 0$}, \label{13} \ee

\nin that holds uniformly for $l\ge 1,$ and making use of the
definition (\ref{sig}) of $\sigma_n,$ proves (\ref{10}).

\nin Observe that uniformity (=the constant implied by the term
$O(\alpha^3l^3)$ in (\ref{13}) does not depend on $l$) is  due  to
the fact that for any real $\alpha,$ $$\vert \frac
{d^s}{d\alpha^s} e^{2\pi i\alpha l}\vert \le (2\pi)^sl^s,\quad
s=1,2,\ldots.$$ \qed

\nin Now we are prepared to deal with the central objective  stated
in the beginning of this section.
 Denote
\be
T=T(n)= \int_0^1 \varphi(\alpha)e^{-2\pi i\alpha n}d\alpha, \quad
n=1,2\ldots \label{11} \ee

\nin  The integrand in (\ref{11}) is  periodic with  period $1.$
So for any $0<\alpha_0\le 1/2,$ the integral $T$ can be written as

\be
T=T_1 + T_2,
\label{12}
\ee

\nin where $T_1=T_1(\alpha_0;n),  T_2=T_2(\alpha_0;n)$ are
integrals of the integrand  in (\ref{11}) over the sets
$[-\alpha_0,+\alpha_0]$ and $[-1/2,-\alpha_0]\cup [\alpha_0,1/2]$
correspondingly. Following the idea of \cite{fr}, \cite{fr1}, we
will first show that for an appropriate choice of $\alpha_0
=\alpha_0(n)$ the main contribution to $T,$ as $n\to \infty$ comes
from $T_1.$ Then, estimating $T_1,$ under $\alpha_0 =\alpha_0(n),
\quad n\to \infty$ we will get the desired asymptotic formula
(\ref{4}).

\nin It is clear  from Lemma \ref{lm3} that the asymptotic
behaviour as $n\to \infty$ of the integral $T_1$ is determined by
the asymptotics of the three key parameters $\sigma_n,$ $B_n^2$
and $\rho_3.$

\nin  First we address the problem for the  class of parameter
functions $a$ of the form
\be
a_j= j^{p-1},\quad  j=1,2,\ldots \quad \rm{for \quad a \quad
given} \quad p> 0. \label{555} \ee It follows from the definition
of the above three parameters that in the case considered the
problem reduces to the estimation of sums of the form
$\sum_{j=1}^n j^k e^{-\sigma j}, \quad k>-1,$ as $n\to \infty.$

\nin To do this we apply
   the integral test for the function
\be
f(u)=(\sigma u)^k e^{-\sigma u}, \quad u\ge0, \quad \sigma\ge 0,
\quad k>-1 .\label{1111} \ee

\nin In the case $k>0$ the function $f$ is strictly increasing  on
$[0,k \sigma^{-1}]$ and is strictly decreasing on
$[k\sigma^{-1},+\infty).$ So, applying the integral test
separately on each of the above intervals we have  in the case
considered
\be
\int _0^{n} f(u)du + \gamma_1f(k\sigma^{-1})\le \sum_{j=1}^n
f(j)\le \int _1^{n}f(u)du + \gamma_2f(k\sigma^{-1}), \label{2022}
\ee where the constants $\gamma_1, \gamma_2$ depend on $k$ only
and $f(k\sigma^{-1})=k^ke^{-k}$.

\nin If now $-1< k \le 0,$ then the function $f$ is strictly
decreasing on $[0;+\infty)$ and the integral test gives

\be
\int _1^{n+1} f(u)du \le \sum_{j=1}^n f(j)\le \int _0^{n}f(u)du .
\label{200} \ee

\nin  Note also that for any $\sigma, b>0$
\be
\int _b^{n}f(u)du =\sigma^k \int _b^{n}u^k e^{-\sigma u} du=
\sigma^{-1} \int _{b\sigma}^{n\sigma }z^k e^{-z}dz .\label {201}
\ee

\nin Then, combining (\ref{1111})-(\ref{201})
we arrive at the desired asymptotic estimate:

$$ \sum_{j=1}^n j^k e^{-\sigma j}\sim \frac{\gamma}{\sigma^{k+1}},
\ \ \rm{as}\ \ 0<\sigma \to 0, \quad n\to \infty,$$

\be \mbox{and}\quad \liminf_{n\to\infty}n\sigma >0 . \label{14}
\ee

\nin In particular, if $n\sigma\to \infty,$ then the constant
$\gamma$ in (\ref{14}) can be found explicitly:
\be
\sum_{j=1}^n j^k e^{-\sigma j}\sim \frac{\Gamma(k+1)}{\sigma^{k+1}},
\label{1499}
\ee
where $\Gamma$ is the gamma function.

\nin Further we will write $\bullet(n)\asymp n^\alpha$ if there
exist positive constants $\gamma_1,\gamma_2, $ s.t.
 $\gamma_1 n^\alpha\le\bullet (n)\le \gamma_2n^\alpha,$ for all
sufficiently large $n.$

\nin Extending (\ref{555}) we consider now the class of functions $a$
satisfying
\be
a_j\asymp j^{p-1},\quad \rm{as}\quad j\to\infty\quad \rm{for \quad
a \quad given} \quad p> 0. \label{655} \ee \nin An obvious
variation of the preceding argument gives in this case  the
following analog of (\ref{14}): $$ 
\sum_{j=1}^n a_j j^l e^{-\sigma j}\asymp
\frac{1}{\sigma^{p+l}},\quad {as} \quad 0<\sigma \to 0, \ \ \ n\to
\infty$$ \be \mbox{and} \quad  \liminf_{n\to\infty}n\sigma>0 ,
\quad l\ge 0 .\label{914} \ee

\nin This immediately implies

\begin{lemma} \hfill
\label{lm4}

\nin Let the function $a$ obey (\ref{655}). Then, as $n\to \infty$,
\be
\sigma_n\asymp n^{-\frac{1}{p+1}}
\label{15}
\ee
\be
B^2_n\asymp n^{\frac{p+2}{p+1}} \label{1555} \ee

\nin and
\be
\rho_3\asymp n^{\frac{p+3}{p+1}}. \ee

\end{lemma}

\nin {\bf Proof:}

\nin By the definition of $\sigma_n,$ Lemma 3 and (\ref{914}), in
this case

\be
n\asymp\frac{1}{\sigma_n^{p+1}}, \quad  \mbox{as \quad $n\to \infty$}
\ee

\nin and, consequently,

\be
\sigma _n\asymp  n^{-\frac{1}{p+1}}, \quad \mbox{as \quad $n\to
\infty$}. \ee

 \nin Now the last two assertions follow from
$(\ref{914})$. \qed

\nin At this point we are prepared to estimate the integral
\be
T_1=T_1(\alpha_0;n)=
\int_{-\alpha_0}^{\alpha_0}\varphi(\alpha)e^{-2\pi i\alpha
n}d\alpha .\ee

\begin{lemma} \hfill
\label{lm5}

\nin Let the function $a$ obey (\ref{655})  and
$\alpha_0=\sigma_n^{\frac{p+2}{2}}\log\ n,\quad p>0.$

\nin Then

\be
T_1(\alpha_0;n)\sim (2\pi B_n^2)^{-1/2}, \quad \mbox{as \quad
$n\to \infty$} .\label{100} \ee
\end{lemma}

\nin {\bf Proof:}

\nin By  (\ref{15}),  $\alpha_0\to 0,$ as $n\to \infty.$
 So, making use of  (\ref{10}) we obtain
\be
T_1 = \int_{-\alpha_0}^{\alpha_0} \exp\left(-2\pi^2\alpha^2B_n^2 +
 O(\alpha^3\rho_3)\right)d\alpha, \quad \mbox{as \quad $n\to
 \infty$}.
\label{19}
\ee

\nin Further, under the condition (\ref{655}), it follows from Lemma 5,
that the above  choice of $\alpha _0$ provides
the following two basic relationships between the parameters
$B^2_n$ and $\rho_3:$
\be
\lim _{n\to\infty} \alpha_0^2 B_n^2=\lim_{n\to \infty} \log^2 n
=+\infty \label{16} \ee \nin and
\be
\lim_{n\to\infty} \alpha_0^3 \rho_3 =\lim_{n\to\infty}
\sigma_n^{\frac{p}{2}}\log^3 n=0 .\label{17} \ee

\nin Since $\lim_{n\to\infty} \alpha^3 \rho_3=0$ for all
$\alpha\in [-\alpha_0,\alpha_0]$, (\ref{19}) and (\ref{16})  imply

\ber T_1\sim \int_{-\alpha_0}^{\alpha_0}
\exp\left(-2\pi^2\alpha^2B_n^2\right)d\alpha=\nonumber\\
\frac{1}{2\pi B_n} \int_{-2\pi\alpha_0B_n}^{2\pi \alpha_0B_n}
\exp(-\frac{z^2}{2}) dz \sim \frac{1}{\sqrt{2\pi B_n^2}}, \quad
{as} \quad n\to \infty .\label{20} \eer \qed

\nin Taking $\alpha_0$ as in Lemma 6,  we write \be
T_2=T_{2,1}+T_{2,2}, \label{T}\ee where
\be
 T_{2,1}= T_{2,1}(\alpha_0;n):=
\int_{ \alpha_0}^{ 1/2} \varphi(\alpha)e^{-2\pi i\alpha n}d\alpha
 \label{eq21} \ee

\nin and  $T_{2,2}$ is the integral of the same integrand, but
over the set $[-1/2,\alpha_0].$

\nin In view of (\ref{T}) and the fact that
$\varphi(-\alpha)=\varphi(\alpha), \ \alpha\in R,$  the rest of
this section is devoted to estimation of the integral $ T_{2,1},$
as $n\to \infty.$

\nin Our starting argument will be the same as in
\cite{fr}. It follows from (\ref{9}) that

\be
\vert \varphi (\alpha)\vert=\exp\left(-2\sum_{j=1}^n
a_je^{-j\sigma_n}\sin^2\pi \alpha j\right), \quad \alpha\in R.
\label{eq22} \ee

\nin Denote by $[x]$ and  $\{x\}$ respectively the integer and
fractional parts
 of a real number $x$ and $\Vert x\Vert$ the distance
from $x$ to the nearest
 integer, so that

\be
\Vert x\Vert=
\left\{
\begin{array}{cc}
\{x\},\quad {if} \quad \{x\}\le 1/2 \\
1-  \{x\},\quad
{if }\quad \{x\}>1/2.
\end{array}
\right.
\label{eq23}
\ee

\nin We will make use of the inequality

\be
\sin\theta \ge \frac{2}{\pi}\theta, \quad 0\le
\theta\le\frac{\pi}{2}. \label{eq24} \ee

\nin Since $\sin^2\pi x= \sin^2\pi \Vert x\Vert$
for any real $x,$ it follows from (\ref{eq24}), (\ref{eq23}) that
for all real $x$

\be
\sin^2\pi x\ge 4 \Vert x\Vert^2.
\ee

\nin Hence,   in view of (\ref{eq22}) and (\ref{eq21})
we have to  estimate the sum

\be
V_n(\alpha):=\sum_{j=1}^n a_je^{-j\sigma_n}\Vert \alpha j\Vert^2,
\quad \alpha_0\le \alpha\le 1/2. \label{eq25} \ee

\begin{lemma} \hfill
\label{lm6}

\nin  Let the function $a$ obey (\ref{655}). Then
\be
V_n(\alpha)\ge \gamma\  \log^2 n, \quad \alpha\in
[\alpha_0,1/2],\quad \mbox{as \quad $n\to \infty$}, \ee \nin where
$\gamma>0.$
\end{lemma}

\nin {\bf Proof:} \nin Take
$\alpha_1=\alpha_1(n)=\sigma_n>\alpha_0(n),$ as $n\to \infty$ and
split the interval $[\alpha_0;1/2]$ into two disjoint subintervals
$I_1:=[\alpha_0,\alpha_1]$ and $I_2:=(\alpha_1,1/2].$ We plan to
prove the assertion separately
 for $\alpha\in I_1$ and $\alpha\in I_2.$


\nin {\bf Interval $I_1.$}

\nin It is clear from the definition (\ref{eq23}) that

\be
\Vert \alpha j\Vert=\alpha j, \quad j\le \frac{1}{2\alpha_1},
\quad \alpha\in I_1. \ee \nin In view of the fact that, by Lemma
3, $(2\alpha_1)^{-1}\to \infty,$ as $n\to \infty,$  while
$(\alpha_1)^{-1}\sigma_n=1,\quad n=1,2, \ldots,$ we apply
(\ref{914}) with $l=2$  to obtain

\be
V_n(\alpha)\ge\alpha^2_0 \sum_{1\le j\le (2\alpha_1)^{-1}} a_j
j^2e^{-j\sigma_n}\asymp \frac{\alpha^2_0}{\sigma_n^{p+2}} =
\log^2 n, \quad \alpha\in I_1. \label{43} \ee

\nin {\bf Interval $I_2.$}

\nin For  a given integer $n$ and a given  $\alpha \in I_2$
 define the
set of integers $$Q(\alpha)=Q(\alpha;n)= \{1\le j \le n:\Vert
\alpha j \Vert \ge 1/4\}.$$ \nin It is clear that
$$Q(\alpha)=\{1\le j\le n: k+1/4\le \alpha j\le k+3/4,\ \
k=0,1,\ldots,\} =$$
 \be\cup_{k=0}^{[\frac{4\alpha
 n-3}{4}]} Q_{k}(\alpha),
\ee \nin  where $Q_k(\alpha)$ denotes the set of integers $\{ j:
\frac{4k+1}{4\alpha}\le j\le \frac{4k+3}{4\alpha}\}.$

\nin Observe that for any $\alpha\in I_2$ and $k\ge 0$
 the set $Q_k(\alpha)$ is not empty, since in this case
 $ \frac{4k+3}{4\alpha}- \frac{4k+1}{4\alpha}\ge 1.$

\nin This yields the following estimate of the sum
$V_n(\alpha),\quad \alpha\in I_2:$

\be
V_n(\alpha)\ge 1/16\sum_{j\in Q(\alpha)}a_je^{-j\sigma_n}= 1/16
\sum_{k=0}^{[\frac{4\alpha n-3}{4}]} \sum_{{j\in Q_k(\alpha)}}
a_je^{-j\sigma_n}. \label{30} \ee

\nin  We now assume that the asymptotic inequality (\ref{655})
holds for all $j\ge N.$ This means that (\ref{655}) is valid for
all $j\in Q_k(\alpha)$ whenever $k\ge \max\{0;\frac{4\alpha
N-1}{4}\}=:K(\alpha). $

\nin We agree, with an obvious abuse of notation, that for a real
$u,$ $Q_u(\alpha)$ is  the interval $$[\frac{4u+1}{4\alpha},
\frac{4u+3}{4\alpha}].$$

\nin Observe that for all sufficiently large $n,$ we have
$0<\alpha^{-1}\sigma_n\le 1,$

\nin $\alpha\in I_2, $ while $n\sigma_n\to\infty$.

\nin Applying now  the integral test to the double sum in the
RHS of (\ref{30}) gives

$$ V_n(\alpha)\ge \gamma \int _{K(\alpha)}^{\frac{4\alpha
n-3}{4}}du \int_{v\in Q_u(\alpha)}v^{p-1}e^{-v\sigma_n}dv=$$
$$\gamma\frac{1}{\sigma_n^p}\int _{K(\alpha)}^{\frac{4\alpha
n-3}{4}}du \int_{v\in (\sigma_n Q_u(\alpha))}v^{p-1}e^{-v}dv\ge$$
$$ \gamma\frac{1}{\sigma_n^p}\int _{K_1(\alpha)}^
{n\sigma_n-\frac{\sigma_n}{2\alpha}}v^{p-1}e^{-v}dv
\int_{\frac{\alpha v}{\sigma_n}-3/4}^{ \frac{\alpha
v}{\sigma_n}-1/4}du \asymp \frac{1}{\sigma_n^p},$$\be \quad
p>0,\quad \mbox{as} \quad n\to \infty, \label{444} \ee

\nin where we denoted
$K_1(\alpha)=\frac{4K(\alpha)+3}{4\alpha}\sigma_n$ and $$\sigma_n
(Q_u(\alpha))=[\sigma_n\frac{4u+1}{4\alpha},
\sigma_n\frac{4u+3}{4\alpha}].$$ \nin Note that the last
inequality in (\ref{444}) is obtained via the change of the order
of integration. Finally, (\ref{444}),  (\ref{15}) and (\ref{43})
prove the claim. \qed \vskip .5cm \nin The last statement of this
section is the desired local limit theorem.

\nin \begin{thm} \hfill

\nin Let the function $a$ obey (\ref{655}). Then
\be
T=Pr(Y=n)\sim \frac{1}{\sqrt{2\pi B_n^2}}, \quad \mbox{as $n\to \infty.$}
\label{47}
\ee
\label{thm123}
\end{thm}

\nin {\bf Proof:}

\nin By Lemma 7 and (\ref{eq21}), (\ref{eq22}), (\ref{T}), we have
\be
T_2\le   e^{- 2\gamma \log^2 n}, \quad \mbox{as $n\to \infty$,}
\label{70} \ee

\nin where $\gamma >0$. In view of (\ref{20}), and (\ref{12}) this
proves (\ref{47}). \qed \vskip .5cm \nin We provide now an
extension of the field of validity of the above local limit
theorem.

\nin We agree to write $ n^{\beta_1}\preceq \bullet(n)\preceq
n^{\beta_2}, \quad 0\le\beta_1\le \beta_2,$ if there exist
positive constants $\gamma_1,\gamma_2, $ s.t. $\gamma_1
n^{\beta_1}\le\bullet (n)\le \gamma_2 n^{\beta_2},$ for all
sufficiently large $n.$

\nin For given $0<p_1\le p_2$ define  the set  ${\cal F}(p_1,
p_2)$ of parameter functions $ a=a(j),\quad  j\in R^+,$ obeying
(\ref{5}) and the condition
\be
j^{p_1-1}\preceq a(j)\preceq j^{p_2-1}, \quad 0< p_1\le p_2. \ee

\nin {\bf Corollary 1}

\nin For an arbitrary $p>0$ and  $0<\epsilon\le \frac{p}{3}$ the
local limit theorem (\ref{47}) is valid for all parameter
functions $a\in {\cal F}(\frac{2p}{3}+\epsilon ; p  ). $

\nin {\bf Proof:}

\nin It is clear from the preceding results that for all
$a\in {\cal F}(p_1;p_2)$, we must have, as $n\to\infty$,

\ber n^{-\frac{1}{p_1+1}}\preceq \sigma_n \preceq
n^{-\frac{1}{p_2+1}}\\ \nonumber \sigma_n^{-(p_1+2)}\preceq B_n^2
\preceq \sigma_n^{-(p_2+2)}\\ \nonumber \sigma_n^{-(p_1+3)}\preceq
\rho_3 \preceq \sigma_n^{-(p_2+3)}. \label{000} \eer

\nin Therefore, setting, as in Lemma 6, $\alpha_0=(B_n)^{-1}\log\
n,$ gives

$$ \alpha_0^3 \rho_3\le\gamma_3 (\log^3
n){\sigma_n}^{\frac{3(p_1+2)}{2}}
\sigma_n^{-(p_2+3)}=$$\be\gamma_3(\log^3 n)
\sigma_n^{\frac{3p_1-2p_2}{2}}\to 0,\quad  \mbox{as}\quad n\to
\infty, \ee \nin provided $3p_1-2 p_2 >0.$

\nin Thus, in the case $a\in {\cal F}(p_1;p_2),$ where $p_1,p_2:
\frac{3p_1}{2}>p_2\ge p_1>0,$ Lemma 6 is valid. The proof of Lemma
7 for this case goes along the same lines, with an obvious
replacement of (\ref{444}) by

\be
\frac{1}{\sigma_n^{p_1}}\preceq V_n(\alpha)\preceq
 \frac{1}{\sigma_n^{p_2}},
\quad \alpha\in I_2. \ee \nin Combining these results  proves the
validity of the local limit theorem  for the  class of functions
$a$ in our  statement. \qed

\nin For our subsequent study we will need the following extension of
(\ref{47}).

\nin {\bf Corollary 2}

\nin Under the conditions of Corollary 1,
\be
Pr(Y=n+h)\sim Pr(Y=n)\sim \frac{1}{\sqrt{2\pi B_n^2}},
\quad \mbox{as $n\to \infty$},
\ee
\nin for a  fixed real $h.$

\nin {\bf Proof:} By (\ref{333}),
\be
\bar T(h;n):=Pr(Y=n+h)= \int_{-1/2}^{1/2} \varphi(\alpha)e^{-2\pi
i\alpha (n+h)}d\alpha, \label{60} \ee where the characteristic
function $\varphi(\alpha)$ is given by (\ref{v1}) and (\ref{v2}).
By Lemma 4, we get

\ber \varphi(\alpha)e^{-2\pi i\alpha (n+h)}&=&
\exp\left(-2\pi^2\alpha^2B_n^2 -2\pi i\alpha h
+O(\alpha^3\rho_3)\right), \nonumber\\ {as}\quad \alpha\to 0,
\label{61} \eer \nin where $\rho_3$ is defined as in (\ref{10}).
Next, let $\alpha_0=\alpha_0(n)$ be as in

\nin  Lemma 6. Denote

\be
\bar T_1(h;n)= \int_{-\alpha_0}^{\alpha_0}\varphi(\alpha)e^{-2\pi
i\alpha (n+h)}d\alpha. \ee

\nin Lemma 6 and  (\ref{61}) imply for a fixed $h\in R$

\be
\bar T_1(h;n) \sim \frac{1}{\sqrt{2\pi B_n^2}}, \quad \mbox{as
$n\to \infty$}. \label{62} \ee

\nin Now it is left to estimate the integral

\be
\bar T_{2}(h;n):= \int_{\alpha_0}^{1/2}\varphi(\alpha)e^{-2\pi
i\alpha (n+h)}d\alpha+
\int_{-1/2}^{-\alpha_0}\varphi(\alpha)e^{-2\pi i\alpha
(n+h)}d\alpha. \ee

\nin Since the function $\varphi(\alpha)$ here is the same as in
(\ref{eq22}), the estimate (\ref{70}) is valid also for $\bar
T_2(h;n),$ which together with (\ref{62}) proves
 the statement.
\qed

\section{The asymptotic formula for $c_n$}
By virtue of  (\ref{3}) and Corollary 1, we obtain the following
asymptotic formula for  $c_n$ valid for all parameter functions
\nin $ a\in {\cal F}(\frac{2p}{3}+\epsilon
; p)  , \quad p>0,\quad \epsilon >0: $
\be
c_n\sim \frac{1}{\sqrt{ 2\pi B_n^2}}
\exp\left(n\sigma_n+\sum_{j=1}^n a_je^{-j\sigma_n}\right),
 \quad \mbox{as $n\to \infty$}.
\label{51}
\ee

\nin {\bf Example.} \nin Let $a_j\asymp j^{p-1}, \quad p> 0,\quad
\mbox{as $j\to \infty$,} .$ Then by Lemma 5, (\ref{914}) and
(\ref{51}),

\be
\log\ c_n\asymp 2n^{\frac{p}{p+1}}-\frac{1}{2} \log 2\pi -
\frac{p+2}{2p+2} \log\ n, \quad \mbox{as $n\to \infty$.}
\label{660} \ee \nin In particular, if $a_j= j^{p-1}, \quad p>0,
\quad j=1,2,\ldots,$ then using (\ref{1499}), the constants in
(\ref{15}), (\ref{1555}) can be found explicitly and we obtain, as
$n\to \infty,$
 \be
\sigma_n\sim\left(\frac{n}{\Gamma(p+1)}\right)^{-\frac{1}{p+1}},\ee
\be B_n^2\sim
\left(\frac{n}{\Gamma(p+1)}\right)^{\frac{p+2}{p+1}}\Gamma(p+2),\ee

\be n\sigma_n\sim \left(\Gamma(p+1)\right)^{\frac{1}{p+1}},\ee

\be \sum_{j=1}^n a_je^{-j\sigma_n}\sim
\Gamma(p)\Big(\frac{n}{\Gamma(p+1)}\Big)^{\frac{1}{p+1}}=
p^{-1}\left(\Gamma(p+1)\right)^{\frac{1}{p+1}}n^{\frac{p}{p+1}} .
 \ee

\nin Hence, (\ref{51}) gives, as $n\to \infty$, \begin{eqnarray}
\log c_n\sim A(p)n^{\frac{p}{p+1}}- \frac{1}{2}\log 2\pi&-&
\frac{p+2}{2p+2} \left(\log \ n - \log\ \Gamma (p+1)
\right)-\nonumber\\1/2 \log\ \Gamma(p+2),  \label{xy}
\end{eqnarray}

\nin where $$A(p)=
(1+p^{-1})\left(\Gamma(p+1)\right)^{\frac{1}{p+1}}.$$ \vskip.5cm

\nin {\bf Remark.} For the two  cases $a_j=$ const and
$a_j=j,\quad j=1,2,\ldots$ the first (=
 the principal ) term in the asymptotic formula
 (\ref{xy})  was obtained in \cite{dgg},
 by solving for large $n$ the corresponding difference equations
(\ref{76}) below. Note that this approach  is not applicable
 even for the class of parameter functions
 $a_j=j^{p-1}, \ p>0$ . \vskip .5cm \nin With the help of
(\ref{51}) we are able to address the question on the validity of
the Conjecture stated in the Section 1 (see
(\ref{rs}),(\ref{rs1})).

\nin {\bf Assertion.}
\nin The conjecture is valid for all parameter functions

\nin $a\in {\cal F}(\frac{2p}{3}+\epsilon ; p),
\quad p>0,\quad \epsilon>0.$

\nin {\bf Proof:}

\nin It is clear from (\ref{1}) that the expression for $c_{n+1}$
can be written
in the following way
\be
c_{n+1}=a_{n+1}+ \bar c_{n+1},
\label{53}
\ee
 where $\bar c_{n+1}$ can be viewed as the value of $c_{n+1}$, (given by
(\ref{1}))
when $a_{n+1}=0.$
 Next recall that the representation (\ref{3}) is true for any real
$\sigma$ and  all $n=1,2,\ldots.$
 We now apply this representation
for the particular case $a_{n+1}=0.$ By (\ref{p}) and (\ref{pr}),
in this case $Pr(X_{n+1}=0)=1$. So, combining (\ref{3}),
(\ref{333}) and taking $\sigma=\sigma_n$  gives

\be
\bar c_{n+1}=\exp\left((n+1)\sigma_n + \sum_{j=1}^n
a_je^{-j\sigma_n}\right)Pr(Y=n+1), \label{vz} \ee

\nin where the random variable $Y$ is defined as in Section 2.

\nin Substituting the expression (\ref{vz}) in (\ref{53}) and
applying (\ref{51}) and Corollary 2, we obtain

\be
 \frac{c_{n+1}}{c_n} \sim e^{\sigma_n} + \frac{a_{n+1}}{c_n},
 \quad \mbox{as $n\to \infty$.}
\label {56}
\ee

In view of the assumption (\ref{5}) and Lemma 3, to complete the
proof we have to show that
\be
\lim_{n\to \infty} \frac {a_{n}}{c_n}=0. \ee

\nin To get this result we will implement the difference equation
(see for references \cite{dgg}, p.460)  derived from (\ref{gen2}):

\ber c_0=1, \quad c_1=a_1, \qquad \nonumber \\ (n+1)c_{n+1} = \sum
_{j = 0} ^{n} (j + 1)a_{j+1}c_{n - j}, \quad n= 1, 2,\ldots
\nonumber \\ \label{76} \eer

\nin This gives for a fixed $ k\ge 1,$
\be
\frac{c_{n+1}}{a_{n+1}}\ge \sum _{j = n-k} ^{n} \frac{(j +
1)a_{j+1}c_{n - j}}{(n + 1)a_{n+1}}. \ee

\nin Consequently, by (\ref{5}) and (\ref{53})

\be
\liminf \frac{c_{n+1}}{a_{n+1}}\ge c_0 +c_1+\ldots+ c_k\ge 1
+a_1+\ldots+ a_k, \ee for any fixed k.

\nin The desired conclusion now follows from the fact
that if

\nin $a\in {\cal F}(\frac{2p}{3}+\epsilon ; p)  , \quad p>0,\quad
\epsilon>0,$ then
\be
\sum _{j = 1} ^{\infty} a_j=\infty. \ee \qed

\nin{\bf Concluding remark}

\nin  Our asymptotic formula (\ref{51}) is restricted to the case

\nin $a\in {\cal F}(\frac{2p}{3}+\epsilon ; p), \quad p>0,\quad
\epsilon>0.$ The reason we require $p>0$ comes from  the fact that
the local limit theorem we proved assumes convergence to the
normal law only. This type of convergence is guaranteed if
(\ref{16}),(\ref{17}) hold, or,
   equivalently, if the parameters
$\rho_3$ and $B_n$ obey the condition:
\be
\lim_{n\to \infty}\frac{\rho_3}{B_n^3}=0. \label{700} \ee

\nin Based on the reasoning preceding (\ref{abc}), it is not
difficult to show that the quantity $\rho_3$ has the following
meaning:
\be
\rho_3=\sum_{l=1}^n E(X_l -EX_l)^3.\label{bor}\ee

\nin In effect, denoting by $U_l$ the
Poisson($\lambda_l:=a_le^{-\sigma l}$) random variable, we have
\be
E(X_l -EX_l)^3=l^3E(U_l -EU_l)^3=l^3\lambda_l,  \ee \nin where the
last step is due to the known property of the third central moment
of Poisson distribution (see, e.g. \cite{LH}, p. 33).

 \nin  (\ref{bor})   explains that (\ref{700}) is  Lyapunov's
sufficient condition for the convergence to the normal law in the
central limit theorem. \nin To demonstrate that for $p\le 0$ the
condition
 (\ref{700}) fails, we consider
 the case $a_j=j^{-1}, \quad j=1,2,\ldots$. It is easy to see from
(\ref{sig}) that for this case $\sigma_n=0,\quad n=1,2\ldots$
Consequently, we have, as $n\to \infty,$ $B_n^2\sim n^2$ and
$\rho_3\sim n^3,$ which gives $\lim_{n\to
\infty}\frac{\rho_3}{B_n^3}=1.$

\nin The above discussion suggests that for $p\le 0$ one can
expect convergence to other stable laws. The study of this case
will require a quite different estimation technique.

{\bf Acknowledgements}

\nin We are very grateful to two referees who carefully read the
paper,
 pointed out a number of inaccuracies, and made many constructive remarks to
 improve the exposition. We are especially thankful to a referee who
 pointed out  the connection of the topic of the paper with the
 field of RCS's.

 The research of the second author was supported by the Fund for
 the Promotion of Research at Technion.

\end{document}